\documentclass[11pt]{amsart}
\usepackage[utf8]{inputenc}
\usepackage[T1]{fontenc}

\usepackage{amsmath,amssymb,amsthm,mathtools}
\usepackage{booktabs}
\usepackage{graphicx}
\usepackage[numbers,sort&compress]{natbib}
\usepackage[margin=2.8cm]{geometry}
\usepackage[dvipsnames]{xcolor}
\usepackage[colorlinks=true,linkcolor=blue!50!black,citecolor=blue!50!black,urlcolor=blue!50!black]{hyperref}

\newtheorem{theorem}{Theorem}[section]
\newtheorem{proposition}[theorem]{Proposition}
\newtheorem{lemma}[theorem]{Lemma}
\newtheorem{corollary}[theorem]{Corollary}
\theoremstyle{definition}
\newtheorem{definition}[theorem]{Definition}
\newtheorem{example}[theorem]{Example}
\theoremstyle{remark}
\newtheorem{remark}[theorem]{Remark}

\newcommand{\R}{\mathbb{R}}
\newcommand{\one}{\mathbf{1}}
\newcommand{\T}{\mathcal{T}}
\newcommand{\eps}{\varepsilon}
\DeclareMathOperator{\RV}{RV}

\title[Comparing and measuring intermediate tail dependence]{Comparing and measuring intermediate tail dependence}
\author{Karl Friedrich Siburg}
\address{Fakult\"at f\"ur Mathematik, Technische Universit\"at Dortmund, 44221 Dortmund, Germany}

\begin{document}

\begin{abstract}
The tail dependence function $\Lambda$ of a bivariate copula, and the tail dependence order based on it, cannot distinguish between copulas whose tail dependence coefficient vanishes; for all of them $\Lambda\equiv 0$. Such copulas, among them the Gaussian copulas with correlation $\rho<1$ and the lower tails of all extreme-value copulas other than the comonotonicity copula, are nevertheless ordered by the much finer information contained in their tail expansion $C(u w)\sim u^{\kappa}\ell(u)\,b(w)$, where $\kappa$ is the tail order, $\ell$ is slowly varying and $b$ is the tail order function of Hua and Joe. We introduce an asymptotic tail order $\preceq$ that refines the tail dependence order and show that it is lexicographic in the three quantities $\kappa$, $\ell$ and $b$. As a consequence, every real-valued measure that is monotone with respect to $\preceq$ is constant on all but countably many tail-order classes; meaningful measures of intermediate tail dependence must therefore be defined level-wise, and we construct such measures. Finally we determine where the fine structure resides: if the joint tail is ray-regular on the exponential scale, then $b$ is determined by the tail order and the one-sided derivatives at $\omega=1/2$ of the angular dependence function of Wadsworth and Tawn; it is a (possibly broken) product of powers. If the joint tail is log-concave on the exponential scale, $b$ is squeezed between two such products.
\end{abstract}


\maketitle

\section{Introduction}\label{sec:intro}

Let $(X,Y)$ be a random pair with continuous margins and copula $C$. The classical (lower) tail dependence coefficient $\lambda_L(C)=\lim_{u\downarrow 0}C(u,u)/u$ \citep{Joe2014} measures the joint extremes along the diagonal. More generally, the tail dependence function
\[
\Lambda(w_1,w_2)=\lim_{u\downarrow 0}\frac{C(uw_1,uw_2)}{u},\qquad w\in[0,\infty)^2,
\]
introduced by \citet{SchmidtStadtmueller2006} and \citet{JoeLiNikoloulopoulos2010}, contains the full first-order information on the lower tail; in particular $\lambda_L(C)=\Lambda(1,1)$. Building on it, \citet{SiburgStrothmannWeiss2024} and, independently, \citet{KoikeKatoHofert2023} ordered copulas by $C_1\le_{td}C_2$ if $\Lambda_1\le\Lambda_2$ pointwise, and \citet{SiburgStrothmannWeiss2024} constructed measures of tail dependence that are monotone with respect to this order, such as the maximal and the average tail dependence. The relation of $\le_{td}$ to local stochastic dominance was clarified by \citet{SiburgStrothmann2024}.

All of this is void for copulas without tail dependence in the usual sense. For these $\Lambda\equiv 0$, so that all of them are equivalent under $\le_{td}$. This class is large and practically relevant: it contains the Gaussian copulas with correlation $\rho<1$, the lower tails of all extreme-value copulas except the comonotonicity copula $M$, the Frank copula and many Archimedean families \citep{HuaJoe2011,Li2024}. For asymptotically independent pairs, \citet{LedfordTawn1996,LedfordTawn1997} and \citet{HuaJoe2011} replaced the first-order expansion by a regularly varying one,
\begin{equation}\label{eq:expansion}
C(uw_1,uw_2)\sim u^{\kappa}\,\ell(u)\,b(w_1,w_2),\qquad u\downarrow 0,
\end{equation}
where $\kappa\ge 1$ is the \emph{tail order} (equal to $1/\eta$ for the coefficient of tail dependence $\eta$ of Ledford and Tawn), $\ell$ is slowly varying at $0$, and the \emph{tail order function} $b$ is homogeneous of degree $\kappa$ with $b(1,1)=1$. The case $1<\kappa<2$ is called \emph{intermediate tail dependence}.

The purpose of this paper is to do for the expansion \eqref{eq:expansion} what \citet{SiburgStrothmannWeiss2024} did for $\Lambda$: to define a consistent order, to understand which measures are compatible with it, and to construct such measures. In their review, \citet{HuaJoe2013} note that a smaller tail order implies a form of stochastic dominance in the joint tail, and they list as a question for future research whether multivariate stochastic orders can be adapted to tail forms in order to compare the strength of tail dependence. Comparisons of two copulas via $(\kappa,\lim\ell)$ along the diagonal were studied by \citet{KoikeKatoYoshiba2025}, who construct a measure and a test of tail equivalence. The order of \citet{Li2013} compares right tails of distributions and is applied to multivariate extreme-value distributions. Our order uses the whole expansion \eqref{eq:expansion}, off-diagonal directions included. Directional tail-weighted dependence measures \citep{LiJoe2024} converge to the tail dependence function and thus degenerate for $\kappa>1$, with a rate determined by the tail order.

More precisely, this paper contains the following.
\begin{enumerate}
\item An \emph{asymptotic tail order} $C_1\preceq C_2$ (Definition~\ref{def:order}) is introduced, which is implied by local stochastic dominance. It refines $\le_{td}$, and it coincides with $\le_{td}$ whenever the dominating copula has positive tail dependence coefficient (Proposition~\ref{prop:td}).
\item On copulas admitting \eqref{eq:expansion}, $\preceq$ is lexicographic: first the tail order $\kappa$ decides, then the slowly varying function $\ell$, and only then the tail order function $b$ (Proposition~\ref{thm:levels}). For instance, for $\rho\in(0,1)$ the Gaussian copula is strictly dominated by the lower Gumbel copula with the same tail order, independently of $b$.
\item As a consequence, any real-valued measure that is monotone in $\preceq$ is constant on all but countably many tail-order classes (Theorem~\ref{thm:impossible}); in particular there is no strictly monotone measure. Measures of intermediate tail dependence therefore have to be defined \emph{level-wise}, and we construct such measures (Section~\ref{sec:measures}).
\item We clarify where the information in $b$ resides (Section~\ref{sec:structure}). Using the angular dependence function $\alpha$ of \citet{WadsworthTawn2013}, we show that $b$ is a product of powers if the joint tail is ray-regular and $\alpha$ is differentiable at $1/2$, and a product broken along the diagonal if $\alpha$ has only one-sided derivatives (Theorem~\ref{thm:rayregular}). The differentiable case is essentially contained in \citet[Sections~3.2 and 4.1]{WadsworthTawn2013}. Under log-concavity of the joint tail on the exponential scale, $b$ is squeezed between two such products without any ray-regularity assumption (Theorem~\ref{thm:logconcave}). For ray-regular tails, and for log-concave tails with $\alpha$ differentiable at $1/2$, the order on the finest level reduces to three real parameters (Corollary~\ref{cor:regularorder}). Broken products occur for the Marshall--Olkin and the generalized Marshall--Olkin copulas; tail order functions beyond (broken) products occur for mixtures and hidden regular variation.
\end{enumerate}

\emph{What is known and what is new.} The tail order, the tail order function and their values for the standard families are due to \citet{LedfordTawn1996,LedfordTawn1997} and \citet{HuaJoe2011}, with further examples in \citet{LiJoe2023} and \citet{Li2024}. The fact that $b$ is a product of powers under ray independence and a differentiable angular dependence function is, in substance, contained in \citet[Sections~3.2 and 4.1]{WadsworthTawn2013}; we only translate it into copula language (Theorem~\ref{thm:rayregular}, differentiable case). New in this paper are the asymptotic tail order and its lexicographic structure (Section~\ref{sec:order}), the impossibility result for monotone measures and the level measures (Section~\ref{sec:measures}), the kinked case of Theorem~\ref{thm:rayregular}, the log-concave sandwich (Theorem~\ref{thm:logconcave}), the reduction of the order to $(\theta,\alpha'_+,\alpha'_-)$ (Corollary~\ref{cor:regularorder}), and the analysis of the generalized Marshall--Olkin and ACIG examples (Example~\ref{ex:families}(4) and Appendix~\ref{app:proofs}).

Throughout we restrict ourselves to the bivariate lower tail; upper tails are treated via the survival copula. We do not discuss estimation in this paper.

\emph{Notation.} We write $\lambda_L(C)$ for the tail dependence coefficient, $\Lambda$ for the tail dependence function, $\Upsilon_L(C)$ for the tail order parameter of \citet{HuaJoe2011} (see Definition~\ref{def:T}; for tail order $1$ it coincides with the tail dependence coefficient, see Lemma~\ref{lem:basic}(4)), and $\alpha$ for the angular dependence function of \citet{WadsworthTawn2013} (see Section~\ref{sec:structure}). For $f,g>0$ we write $f\sim g$ if $f/g\to 1$. $\RV_0(0^+)$ is the set of functions slowly varying at $0$, and $\one=(1,1)$.

\emph{Use of AI tools.} The AI model Claude (Anthropic) was used for \LaTeX\ editing, generating Figure 1, and for checking the calculations in the Appendix. The author has checked all statements, proofs and references
and takes full responsibility for the content.

\section{Tail orders and tail order functions}\label{sec:prelim}

\begin{definition}\label{def:T}
Let $\T$ denote the set of bivariate copulas $C$ for which there exist $\kappa=\kappa_C>0$ and $\ell=\ell_C\in\RV_0(0^+)$ such that
\[
b_C(w):=\lim_{u\downarrow 0}\frac{C(uw_1,uw_2)}{u^{\kappa}\ell(u)}
\]
exists in $[0,\infty)$ for all $w\in(0,\infty)^2$, with $b_C(\one)=1$. We call $\kappa_C$ the \emph{tail order}, $\ell_C$ the \emph{slowly varying part} and $b_C$ the \emph{tail order function} of $C$ \citep{HuaJoe2011}. If $\Upsilon_L(C):=\lim_{u\downarrow0}\ell_C(u)\in[0,\infty]$ exists, it is the \emph{tail order parameter}.
\end{definition}

The tail order is unique, and $\ell_C$ is unique up to asymptotic equivalence, because $C(u,u)\sim u^{\kappa}\ell(u)$. The following properties are elementary; see \citet[Proposition~1]{HuaJoe2011}.

\begin{lemma}\label{lem:basic}
Let $C\in\T$ with tail order $\kappa$ and tail order function $b$. Then:
\begin{enumerate}
\item $\kappa\ge1$;
\item $b$ is homogeneous of degree $\kappa$, nondecreasing in each argument and $2$-increasing;
\item $\min(w_1,w_2)^{\kappa}\le b(w)\le \max(w_1,w_2)^{\kappa}$ for all $w\in(0,\infty)^2$; in particular $b>0$ on $(0,\infty)^2$;
\item if $\kappa>1$ then $\Lambda_C\equiv0$; if $\kappa=1$ and $\Upsilon_L(C)$ exists, then $\Lambda_C=\Upsilon_L(C)\,b$, so that $\Upsilon_L(C)=\Lambda_C(\one)=\lambda_L(C)$.
\end{enumerate}
\end{lemma}

\begin{proof}
(1) follows from $C(u,u)\le u$. Homogeneity follows from $\ell(tu)/\ell(u)\to1$; monotonicity and $2$-increasingness are inherited from $C$. For (3) use monotonicity and homogeneity: $b(w)\ge b(\min(w)\one)=\min(w)^{\kappa}$, and similarly for the upper bound. (4) is immediate from $C(uw)/u=u^{\kappa-1}\ell(u)\,b(w)(1+o(1))$ and $u^{\kappa-1}\ell(u)\to0$ for $\kappa>1$.
\end{proof}

For $\kappa>1$ the set of possible tail order functions is large. By \citet[Theorem~8.2.3]{Jaworski2010}, every function on $[0,\infty)^2$ that is homogeneous of degree $\kappa>1$, grounded, $2$-increasing and locally Lipschitz coincides near the origin with some copula. In particular, if $\tau$ is such a function then so is $c\tau$ for every $c>0$: tail order functions multiplied by tail order parameters are not bounded from above. This is a fundamental difference to the case $\kappa=1$, where $\Lambda\le\min$.

\begin{example}\label{ex:families}
We collect the copulas used throughout. All statements refer to the lower tail; when an upper tail is meant, we pass to the survival copula $\widehat C$.
\begin{enumerate}
\item \emph{Extreme-value copulas.} Let $C(u,v)=\exp\{-A(-\log u,-\log v)\}$ with a convex, $1$-homogeneous $A$ satisfying $\max(s)\le A(s)\le s_1+s_2$. Then $C(u,u)=u^{A(\one)}$, so $\kappa=A(\one)$ and $\ell\equiv1$. By convexity, $A(\one)=1$ forces $A=\max$, so $\kappa>1$ (and hence $\Lambda\equiv0$) for every extreme-value copula except $M(u,v)=\min(u,v)$. If $A$ is differentiable at $\one$ then $b(w)=w_1^{\partial_1A(\one)}w_2^{\partial_2A(\one)}$ \citep[Example~2]{HuaJoe2011}. The (lower) Gumbel copula, $A(s)=(s_1^{\delta}+s_2^{\delta})^{1/\delta}$ with $\delta\ge1$, has $\kappa=2^{1/\delta}$ and $b(w)=(w_1w_2)^{\kappa/2}$. The Marshall--Olkin copula $C(u,v)=uv\max(u,v)^{-\beta}$, $\beta\in(0,1)$, corresponds to $A(s)=s_1+s_2-\beta\min(s_1,s_2)$, which is not differentiable at $\one$; here $\kappa=2-\beta$ and
\[
b(w)=w_1w_2\max(w_1,w_2)^{-\beta}=\min\bigl(w_1w_2^{\kappa-1},\,w_1^{\kappa-1}w_2\bigr).
\]
\item \emph{Gaussian copulas.} For $\rho\in(0,1)$ the Gaussian copula has $\kappa=2/(1+\rho)$ and $\ell(u)\propto(-\log u)^{-\rho/(1+\rho)}$ \citep[Example~1]{HuaJoe2011}, and $b(w)=(w_1w_2)^{\kappa/2}$ \citep[Appendix~B]{LiJoe2023}. Lemma~\ref{lem:gauss} in the Appendix gives a self-contained proof for all $\rho\in(-1,1)$ together with the constant:
\begin{equation}\label{eq:gaussell}
\ell(u)\sim c_\rho\,\bigl(2\log(1/u)\bigr)^{-\rho/(1+\rho)},\qquad c_\rho=\frac{(1+\rho)^{3/2}}{(1-\rho)^{1/2}}\,(2\pi)^{-\rho/(1+\rho)}.
\end{equation}
In particular $\ell(u)\to0$ for $\rho>0$.
\item \emph{Archimedean copulas with light-tailed frailty.} Here $C_\psi(u,v)=\psi(\psi^{-1}(u)+\psi^{-1}(v))$, where the generator $\psi(s)=E[e^{-sV}]$ is the Laplace transform of a positive random variable $V$, the frailty. If $\psi(s)\sim a_1s^{q}e^{-a_2s^{r}}$ as $s\to\infty$ with $a_1,a_2>0$ and $0<r<1$ (and $\psi'$ correspondingly), then $\kappa=2^{r}$, $\ell(u)\propto(-\log u)^{\zeta}$ with $\zeta=(q/r)(1-2^{r})$, and $b(w)=(w_1w_2)^{\kappa/2}$ \citep[Proposition~6]{HuaJoe2011}; see also \citet[Theorem~3.3]{Li2024} for rapidly varying generators, which include the Frank copula with $\kappa=2$ and $b(w)=w_1w_2$. Depending on the sign of $q$, $\ell$ tends to $0$, to a constant or to $\infty$.
\item \emph{Generalized Marshall--Olkin copulas.} Let $Z_1,Z_2,R$ be independent with $P(Z_i>z)=1/z$ and $P(R>r)=r^{-\kappa}$ for $z,r\ge1$, where $1<\kappa<2$, and put $X_i=\max(Z_i,R)$. Then $(1/X_1,1/X_2)=(\min(1/Z_1,1/R),\min(1/Z_2,1/R))$ has a generalized Marshall--Olkin distribution in the sense of \citet{LiPellerey2011}, with $1/R$ as the shock common to both components. Let $C^{\mathrm{GMO}}_\kappa$ be the copula of $(1/X_1,1/X_2)$, that is, the survival copula of $(X_1,X_2)$, or equivalently the distribution function of $(\bar F(X_1),\bar F(X_2))$ with $\bar F(x)=P(X_i>x)$. (The survival copula of $(1/X_1,1/X_2)$ is a generalized Marshall--Olkin survival copula in the terminology of \citet{LiPellerey2011}; for exponential shocks these reduce to the Marshall--Olkin copulas of item~(1).) Then $C^{\mathrm{GMO}}_\kappa\in\T$ with tail order $\kappa$, $\ell\to1$ and
\[
b(w)=\min(w_1,w_2)^{\kappa}.
\]
Indeed, $\bar F(x)=x^{-1}+x^{-\kappa}-x^{-\kappa-1}$, so $\bar F^{-1}(u)=u^{-1}(1+O(u^{\kappa-1}))$. For $x_1,x_2\ge1$ with $m=\min(x_1,x_2)$ and $M=\max(x_1,x_2)$,
\[
P(X_1>x_1,X_2>x_2)=M^{-\kappa}+\frac{m^{-\kappa}-M^{-\kappa}}{M}+\frac{1-m^{-\kappa}}{x_1x_2}.
\]
Inserting $x_i=\bar F^{-1}(uw_i)$, the first term equals $u^{\kappa}\min(w)^{\kappa}(1+o(1))$, the second is $O(u^{\kappa+1})$ and the third is $O(u^2)=o(u^{\kappa})$. More generally, $X_i=\max(Z_i,c_iR)$ with $c_1,c_2>0$ gives $b(w)=\min(c_1w_1,c_2w_2)^{\kappa}/\min(c_1,c_2)^{\kappa}$.
\item \emph{Upper tail of the ACIG copula.} The Archimedean copula $C_\psi$ generated by the Laplace transform of an inverse Gamma variable with shape $\kappa\in(1,2)$ has upper tail order $\kappa$ \citep[Example~4]{HuaJoe2011}. Its survival copula $C^{\mathrm{ig}}_\kappa:=\widehat C_\psi$ has
\[
b(w)=\frac{(w_1+w_2)^{\kappa}-w_1^{\kappa}-w_2^{\kappa}}{2^{\kappa}-2}
\]
\citep[Proposition~5]{HuaJoe2011}, and
\[
\ell\to\Upsilon_L(C^{\mathrm{ig}}_\kappa)=c_h(2^{\kappa}-2)(\kappa-1)^{\kappa},\qquad c_h=-\frac{\Gamma(1-\kappa)}{\Gamma(1+\kappa)}>0
\]
(Lemma~\ref{lem:acig}). For $\kappa=3/2$ this gives $c_h=8/3$ and $\Upsilon_L(C^{\mathrm{ig}}_{3/2})=(8-4\sqrt2)/3\approx0.7810$.
\item \emph{Mixtures.} If $C_1,C_2\in\T$ have the same tail order and $\ell_1\sim\ell_2$, then $pC_1+(1-p)C_2\in\T$ with tail order function $pb_1+(1-p)b_2$, for $p\in[0,1]$.
\end{enumerate}
\end{example}

\section{The asymptotic tail order}\label{sec:order}

\begin{definition}\label{def:order}
Let $\mathcal C_+$ be the set of bivariate copulas $C$ with $C>0$ on $(0,1]^2$. For $C_1,C_2\in\mathcal C_+$ we write $C_1\preceq C_2$ (\emph{$C_1$ is asymptotically less lower tail dependent than $C_2$}) if
\[
\limsup_{u\downarrow0}\frac{C_1(uw_1,uw_2)}{C_2(uw_1,uw_2)}\le1\qquad\text{for all }w\in(0,\infty)^2.
\]
We write $C_1\prec C_2$ if $C_1\preceq C_2$ but not $C_2\preceq C_1$, and $C_1\sim_t C_2$ if both $C_1\preceq C_2$ and $C_2\preceq C_1$.
\end{definition}

Restricting to $\mathcal C_+$ is necessary: for the countermonotonicity copula $W(u,v)=\max(u+v-1,0)$ the quotient $W(uw)/W(uw)$ is undefined for small $u$. It is no restriction for our purposes, because $\T\subseteq\mathcal C_+$: for $C\in\T$ we have $C(t,t)>0$ for small $t$ by \eqref{eq:expansion} and Lemma~\ref{lem:basic}(3), hence for all $t\in(0,1]$ by monotonicity, and $C(v)\ge C(\min(v),\min(v))>0$ for $v\in(0,1]^2$.

Recall that $C_1\le_{loc}C_2$ (local stochastic dominance) means $C_1\le C_2$ on $[0,\delta]^2$ for some $\delta>0$ \citep{SiburgStrothmann2024}.

\begin{proposition}\label{prop:basicorder}
\begin{enumerate}
\item $\preceq$ is a preorder on $\mathcal C_+$, i.e.\ reflexive and transitive.
\item $C_1\le_{loc}C_2$ implies $C_1\preceq C_2$; in particular the concordance order implies $\preceq$.
\item $\preceq$ is neither antisymmetric nor total, even on $\T$.
\end{enumerate}
\end{proposition}

\begin{proof}
(1) Reflexivity is clear. If $C_1\preceq C_2\preceq C_3$, then $\limsup C_1/C_3\le\limsup(C_1/C_2)\cdot\limsup(C_2/C_3)\le1$, since both factors have finite $\limsup$. (2) is obvious. (3) Any two copulas that agree near the origin are equivalent; non-totality follows from Example~\ref{ex:chain} below.
\end{proof}

Recall that $C_1\le_{td}C_2$ (tail dependence order) means $\Lambda_1\le\Lambda_2$ on $[0,\infty)^2$, where $\Lambda_i$ is the tail dependence function of $C_i$ \citep[Definition~2.1]{SiburgStrothmannWeiss2024}; this order is defined for copulas admitting a tail dependence function.

\begin{proposition}\label{prop:td}
Let $C_1,C_2\in\mathcal C_+$ admit tail dependence functions $\Lambda_1,\Lambda_2$.
\begin{enumerate}
\item If $C_1\preceq C_2$, then $C_1\le_{td}C_2$.
\item If $\lambda_L(C_2)=\Lambda_2(\one)>0$, i.e.\ if $C_2$ has a positive tail dependence coefficient, then $C_1\le_{td}C_2$ implies $C_1\preceq C_2$.
\end{enumerate}
Hence $\preceq$ refines $\le_{td}$, and both orders coincide as far as comparisons with copulas of positive tail dependence coefficient are concerned.
\end{proposition}

\begin{proof}
(1) Write $C_1(uw)/u=\bigl(C_1(uw)/C_2(uw)\bigr)\bigl(C_2(uw)/u\bigr)$. The first factor has $\limsup\le1$ and the second converges to $\Lambda_2(w)$, so $\Lambda_1(w)\le\Lambda_2(w)$. (2) Monotonicity and homogeneity of $\Lambda_2$ give $\Lambda_2(w)\ge\lambda_L(C_2)\min(w)>0$ for all $w\in(0,\infty)^2$, so $C_1(uw)/C_2(uw)\to\Lambda_1(w)/\Lambda_2(w)\le1$.
\end{proof}

The next result is the basic structural fact: on $\T$ the order $\preceq$ is lexicographic in tail order, slowly varying part and tail order function.

\begin{proposition}\label{thm:levels}
Let $C_1,C_2\in\T$ with tail orders $\kappa_i$, slowly varying parts $\ell_i$ and tail order functions $b_i$, and put $\theta:=\limsup_{u\downarrow0}\ell_1(u)/\ell_2(u)\in[0,\infty]$. Then $C_1\preceq C_2$ if and only if one of the following holds:
\begin{enumerate}
\item $\kappa_1>\kappa_2$;
\item $\kappa_1=\kappa_2$ and $\theta=0$;
\item $\kappa_1=\kappa_2$, $0<\theta<\infty$, and $\theta\,b_1\le b_2$ on $(0,\infty)^2$.
\end{enumerate}
\end{proposition}

\begin{proof}
For fixed $w$,
\[
\frac{C_1(uw)}{C_2(uw)}=u^{\kappa_1-\kappa_2}\,\frac{\ell_1(u)}{\ell_2(u)}\,\frac{b_1(w)}{b_2(w)}\,(1+o(1)),
\]
where $0<b_i(w)<\infty$ by Lemma~\ref{lem:basic}. If $\kappa_1\ne\kappa_2$, then $u^{\kappa_1-\kappa_2}\ell_1(u)/\ell_2(u)$ tends to $0$ if $\kappa_1>\kappa_2$ and to $\infty$ if $\kappa_1<\kappa_2$, because $\ell_1/\ell_2$ is slowly varying \citep[Proposition~1.3.6]{BinghamGoldieTeugels1987}. If $\kappa_1=\kappa_2$, the $\limsup$ of the ratio equals $\theta\,b_1(w)/b_2(w)$. This is $\le1$ for all $w$ if and only if $\theta=0$, or $\theta\in(0,\infty)$ and $\theta b_1\le b_2$.
\end{proof}

\begin{remark}\label{rem:levels}
(a) Proposition~\ref{thm:levels} makes precise the remark of \citet{HuaJoe2011,HuaJoe2013} that a smaller tail order means stronger tail dependence, and the observation of \citet{KoikeKatoYoshiba2025} that tail order parameters are comparable only for equal tail orders. Whenever the measure $\xi_{C_1,C_2,w}$ of \citet{KoikeKatoYoshiba2025} is defined, $C_1\preceq C_2$ implies $\xi_{C_1,C_2,w}\ge0$, because $\xi$ only uses the diagonal.

(b) The middle level is not an artefact. Copulas with the same tail order often have different slowly varying parts, and then $b$ is irrelevant for the comparison (Example~\ref{ex:chain}).

(c) On $\{C\in\T:\kappa_C=1,\ \Upsilon_L(C)>0\}$ the order $\preceq$ is just $\le_{td}$, by Lemma~\ref{lem:basic}(4) and Proposition~\ref{prop:td}.
\end{remark}

\begin{example}\label{ex:chain}
Fix $\kappa=3/2$ and consider the following copulas with tail order $3/2$: the lower Gumbel copula $C^{\mathrm{Gu}}$ ($\delta=\log2/\log(3/2)$), the Marshall--Olkin copula $C^{\mathrm{MO}}$ ($\beta=1/2$), the generalized Marshall--Olkin copula $C^{\mathrm{GMO}}_{3/2}$, the ACIG survival copula $C^{\mathrm{ig}}_{3/2}$ and the Gaussian copula $C^{\mathrm{Ga}}$ with $\rho=1/3$. The first three have $\ell\to1$, so the same tail order parameter $\Upsilon_L(C)=1$, and their tail order functions satisfy, for $r=w_1/w_2\le1$ (by symmetry this suffices),
\[
\frac{b^{\mathrm{GMO}}(w)}{b^{\mathrm{MO}}(w)}=r^{\kappa-1}\le1,\qquad\frac{b^{\mathrm{MO}}(w)}{b^{\mathrm{Gu}}(w)}=r^{1-\kappa/2}\le1 .
\]
The Gaussian copula has $\ell\to0$ by \eqref{eq:gaussell}. Hence, by Proposition~\ref{thm:levels},
\[
C^{\mathrm{Ga}}\prec C^{\mathrm{GMO}}_{3/2}\prec C^{\mathrm{MO}}\prec C^{\mathrm{Gu}},
\]
although $C^{\mathrm{GMO}}$, $C^{\mathrm{MO}}$ and $C^{\mathrm{Gu}}$ share the tail order and the tail order parameter. In contrast, $C^{\mathrm{ig}}_{3/2}$ is not comparable with $C^{\mathrm{GMO}}_{3/2}$ or $C^{\mathrm{MO}}$. Indeed $\Upsilon_L(C^{\mathrm{ig}}_{3/2})=(8-4\sqrt2)/3<1$, but $\Upsilon_L(C^{\mathrm{ig}}_{3/2})\,b^{\mathrm{ig}}(r,1)\sim\Upsilon_L(C^{\mathrm{ig}}_{3/2})\,\kappa r/(2^{\kappa}-2)=\sqrt2\,r$ as $r\to0$, which eventually exceeds both $b^{\mathrm{GMO}}(r,1)=r^{3/2}$ and $b^{\mathrm{MO}}(r,1)=r$. On the other hand $C^{\mathrm{ig}}_{3/2}\prec C^{\mathrm{Gu}}$ by Lemma~\ref{lem:acigineq}(2). Finally, an Archimedean copula as in Example~\ref{ex:families}(3) with $r=\log_2(3/2)$ and $q<0$ has $\ell\to\infty$ and strictly dominates all copulas above.
\end{example}

The asymptotic tail order controls the copulas themselves near the origin, uniformly on cones away from the axes. This is the analogue of the relation between tail dependence order and local stochastic dominance studied by \citet{SiburgStrothmann2024}.

\begin{proposition}\label{prop:cone}
Let $C_1,C_2\in\T$ have the same tail order, $\ell_1/\ell_2\to\theta\in[0,\infty)$, continuous tail order functions, and $\theta b_1<b_2$ on $(0,\infty)^2$. Then for every $\eps\in(0,1)$ there is $\delta>0$ such that
\[
C_1(v)<C_2(v)\qquad\text{for all }v\in(0,\delta)^2\text{ with }\eps\le v_1/v_2\le1/\eps .
\]
\end{proposition}

\begin{proof}
Let $K=\{w:w_1+w_2=1,\ \eps/(1+\eps)\le w_1\le1/(1+\eps)\}$ and $f_u^{(i)}(w)=C_i(uw)/(u^{\kappa}\ell_2(u))$. Each $f^{(i)}_u$ is nondecreasing in both coordinates and converges pointwise to the continuous limit $\theta b_1$ resp.\ $b_2$. Bracketing $w$ between grid points and using uniform continuity of the limit shows that the convergence is uniform on compact subsets of $(0,\infty)^2$, in particular on $K$. Since $\eta:=\min_K(b_2-\theta b_1)>0$, there is $u_0$ with $f^{(1)}_u<f^{(2)}_u$ on $K$ for $u<u_0$. Every $v$ in the cone with $v_1+v_2<u_0$ can be written $v=uw$ with $u=v_1+v_2$ and $w\in K$; take $\delta=u_0/2$.
\end{proof}

\section{Measures of intermediate tail dependence}\label{sec:measures}

Following \citet[Definition~3.1]{SiburgStrothmannWeiss2024}, a \emph{measure} on a subclass $\mathcal S\subseteq\T$ is a map $\mu:\mathcal S\to\R$ with $C_1\preceq C_2\Rightarrow\mu(C_1)\le\mu(C_2)$. For $\kappa\ge1$ let $\T_\kappa=\{C\in\T:\kappa_C=\kappa\}$.

\subsection{An impossibility result}

\begin{theorem}\label{thm:impossible}
Let $\mu:\T\to\R$ be a measure. Then the set $N(\mu)=\{\kappa\ge1:\mu\text{ is not constant on }\T_\kappa\}$ is at most countable.
\end{theorem}

\begin{proof}
For $\kappa\in N(\mu)$ let $J_\kappa=(\inf\mu(\T_\kappa),\sup\mu(\T_\kappa))$, a nonempty open interval. If $\kappa<\kappa'$, then $C'\preceq C$ for all $C'\in\T_{\kappa'}$ and $C\in\T_{\kappa}$ by Proposition~\ref{thm:levels}(1). Hence $\sup\mu(\T_{\kappa'})\le\inf\mu(\T_\kappa)$, so $J_{\kappa'}\cap J_\kappa=\emptyset$. A family of pairwise disjoint nonempty open intervals in $\R$ is at most countable, since each contains a rational number.
\end{proof}

\begin{corollary}\label{cor:nostrict}
There is no strictly monotone measure on $\T$, i.e.\ no $\mu:\T\to\R$ with $C_1\prec C_2\Rightarrow\mu(C_1)<\mu(C_2)$.
\end{corollary}

\begin{proof}
For every $\kappa\in(1,2)$ we have $C^{\mathrm{GMO}}_\kappa\prec C^{\mathrm{Gu}}_\kappa$, where $C^{\mathrm{Gu}}_\kappa$ is the Gumbel copula with $2^{1/\delta}=\kappa$: both have $\ell\to1$, and $\min(w)^{\kappa}\le(w_1w_2)^{\kappa/2}$ with strict inequality off the diagonal. A strictly monotone measure would be nonconstant on uncountably many classes $\T_\kappa$, contradicting Theorem~\ref{thm:impossible}.
\end{proof}

Theorem~\ref{thm:impossible} is a variant of the classical fact that the lexicographic order on $\R^2$ admits no real-valued representation \citep{Debreu1954}. It explains why the problem does not arise for the tail dependence order: there all nontrivial information sits in the single class $\kappa=1$, $\Upsilon_L(C)>0$. For intermediate tail dependence there are two consistent ways out.
\begin{itemize}
\item \emph{Global measures} depend on the tail order only, $\mu(C)=\varphi(\kappa_C)$ with $\varphi$ nonincreasing. An example is $\bar\chi=2/\kappa-1$ of \citet{ColesHeffernanTawn1999}. By Proposition~\ref{thm:levels}, $C_1\preceq C_2$ implies $\kappa_1\ge\kappa_2$, so these are measures.
\item \emph{Level measures} compare copulas within one class of equal tail order and comparable slowly varying parts. They are the subject of the next subsection.
\end{itemize}
A complete numerical description of the tail is then a lexicographically ordered triple consisting of the tail order, the slowly varying part and a level measure.

\subsection{Level measures}

\begin{definition}\label{def:level}
For $\kappa\ge1$ and $L_0\in\RV_0(0^+)$ the \emph{level} $\T(\kappa,L_0)$ is the set of $C\in\T_\kappa$ with $\ell_C/L_0\to\theta_C\in(0,\infty)$. For $C\in\T(\kappa,L_0)$ the \emph{tail profile} is
\[
\tau_C(w)=\lim_{u\downarrow0}\frac{C(uw)}{u^{\kappa}L_0(u)}=\theta_C\,b_C(w),\qquad t_C(s)=\tau_C(s,1-s),\ s\in[0,1].
\]
For $w$ on a coordinate axis the limit is trivial: $C(0,v)=C(v,0)=0$ gives $\tau_C(0,w_2)=\tau_C(w_1,0)=0$, and in particular $t_C(0)=t_C(1)=0$. Thus $\tau_C=\theta_C\,b_C$ on $(0,\infty)^2$, while $\tau_C=0$ on the axes.
\end{definition}

By homogeneity $\tau_C(w)=(w_1+w_2)^{\kappa}t_C(w_1/(w_1+w_2))$, so $t_C$ determines $\tau_C$, just as the tail dependence function restricted to $\{s+t=1\}$ in \citet[Remark~2.3]{SiburgStrothmannWeiss2024}. By Proposition~\ref{thm:levels},
\begin{equation}\label{eq:levelorder}
C_1\preceq C_2\iff\tau_{C_1}\le\tau_{C_2}\iff t_{C_1}\le t_{C_2}\qquad(C_1,C_2\in\T(\kappa,L_0)).
\end{equation}
Replacing $L_0$ by $cL_0$ multiplies all profiles by $1/c$, so rankings by positively homogeneous level measures do not depend on the choice of $L_0$. By Lemma~\ref{lem:basic}(3),
\begin{equation}\label{eq:tbounds}
\theta_C\min(s,1-s)^{\kappa}\le t_C(s)\le\theta_C\max(s,1-s)^{\kappa}\le\theta_C,
\end{equation}
so $t_C$ is bounded and Borel measurable, being a pointwise limit of continuous functions.

\begin{proposition}\label{thm:levelmeasures}
On every level $\T(\kappa,L_0)$ the following maps are measures with values in $[0,\infty)$:
\begin{enumerate}
\item the \emph{normalised tail order parameter} $\tau_C(\one)=2^{\kappa}t_C(1/2)=\theta_C$ (equal to $\Upsilon_L(C)$ if $L_0\equiv1$) and, more generally, $t_C(s_0)$ for fixed $s_0\in(0,1)$;
\item the \emph{maximal intermediate tail dependence} $\mu_\infty(C)=\sup_{s\in[0,1]}t_C(s)$;
\item the \emph{average intermediate tail dependence} $\mu_1(C)=\int_0^1t_C(s)\,ds$;
\item the $L^p$-norms $\mu_p(C)=\bigl(\int_0^1t_C(s)^p\,ds\bigr)^{1/p}$, $1\le p<\infty$.
\end{enumerate}
More generally, $\Phi(t_C)$ is a measure for every functional $\Phi$ that is nondecreasing with respect to the pointwise order.
\end{proposition}

\begin{proof}
All quantities are finite by \eqref{eq:tbounds}, and monotonicity follows from \eqref{eq:levelorder}.
\end{proof}

\begin{remark}\label{rem:differences}
Three differences to the case $\kappa=1$ \citep{SiburgStrothmannWeiss2024} are worth noting.
\begin{enumerate}
\item The level measures are unbounded for $\kappa>1$. If $\tau$ is a locally Lipschitz tail profile, then so is $c\tau$ for every $c>0$ \citep[Theorem~8.2.3]{Jaworski2010}, so no normalisation to $[0,1]$ exists.
\item $t_C$ need not be concave, and for exchangeable $C$ the maximum need not be attained at $s=1/2$, in contrast to \citet[Proposition~3.6]{SiburgStrothmannWeiss2024}. Take the generalized Marshall--Olkin copula of Example~\ref{ex:families}(4) with weights $(c_1,c_2)=(3,1)$ and $\kappa=3/2$ and mix it with its transpose. The result is an exchangeable copula with profile $t(s)=\tfrac12[\min(3s,1-s)^{3/2}+\min(s,3(1-s))^{3/2}]$, for which $t(1/4)=(1+3\sqrt3)/16\approx0.387>t(1/2)=2^{-3/2}\approx0.354$.
\item The tail order parameter is the analogue of the tail dependence coefficient. Like the latter, it may fail to distinguish copulas that are strictly ordered (Example~\ref{ex:chain}).
\end{enumerate}
\end{remark}

\begin{figure}[ht]
\centering
\includegraphics[width=0.85\textwidth]{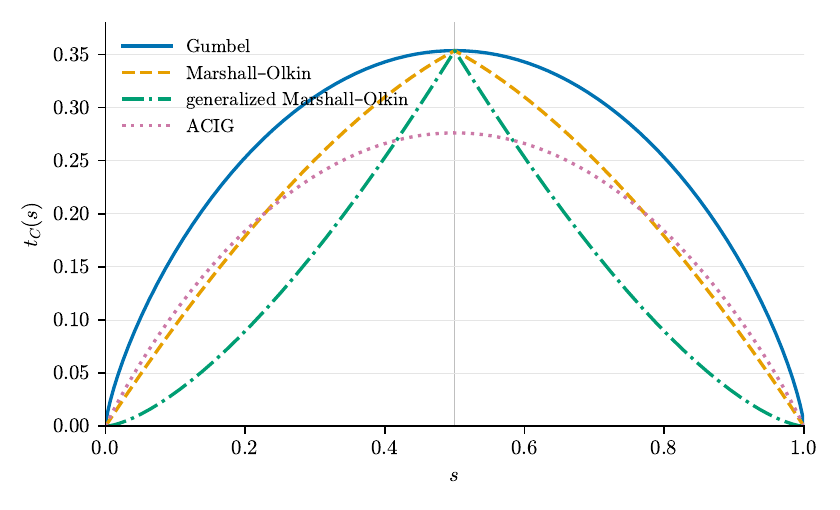}
\caption{Tail profiles $t_C$ on the level $\T(3/2,1)$. The Gumbel, the Marshall--Olkin and the generalized Marshall--Olkin copula coincide at $s=1/2$ (equal tail order parameter) but are strictly ordered, $C^{\mathrm{GMO}}\prec C^{\mathrm{MO}}\prec C^{\mathrm{Gu}}$. The ACIG profile lies below the Gumbel profile but crosses the other two, so it is incomparable with them.}\label{fig:profiles}
\end{figure}

\begin{example}\label{ex:table}
Figure~\ref{fig:profiles} shows the tail profiles of the copulas of Example~\ref{ex:chain}, all on the level $\T(3/2,1)$, and Table~\ref{tab:measures} lists their level measures. The Gaussian copula with $\rho=1/3$ belongs to the lower level $\T\bigl(3/2,(\log(1/u))^{-1/4}\bigr)$. The tail order parameter does not separate the Gumbel, the Marshall--Olkin and the generalized Marshall--Olkin copula, whereas the average does, in accordance with the chain in Example~\ref{ex:chain}. For the incomparable pair $C^{\mathrm{ig}}_{3/2}$ and $C^{\mathrm{GMO}}_{3/2}$ the measures disagree: the maximum ranks $C^{\mathrm{GMO}}_{3/2}$ higher, the average ranks the ACIG copula higher. This is the inconsistency phenomenon of \citet[Fig.~3]{SiburgStrothmannWeiss2024} on the level of intermediate tail dependence.
\end{example}

\begin{table}[ht]
\centering
\caption{Level measures on $\T(3/2,1)$: exact values and decimals (rounded to four digits). In all four cases the maximum is attained at $s=1/2$. For the ACIG copula, $\Upsilon_L(C^{\mathrm{ig}}_{3/2})/(2^{3/2}-2)=2\sqrt2/3$.}\label{tab:measures}
\small
\begin{tabular}{lcccc}
\toprule
copula & $t_C(s)$, $s\le1/2$ & $\mu_\infty=t_C(1/2)$ & $\mu_1$ & $\mu_2^2$\\
\midrule
Gumbel & $(s(1-s))^{3/4}$ & $2^{-3/2}\approx0.3536$ & $\frac{3\Gamma(3/4)^2}{10\sqrt\pi}\approx0.2542$ & $\frac{3\pi}{128}\approx0.0736$\\[3pt]
Marshall--Olkin & $s(1-s)^{1/2}$ & $2^{-3/2}\approx0.3536$ & $\frac{8}{15}-\frac{7\sqrt2}{30}\approx0.2034$ & $\frac{5}{96}\approx0.0521$\\[3pt]
gen.\ Marshall--Olkin & $s^{3/2}$ & $2^{-3/2}\approx0.3536$ & $\frac{\sqrt2}{10}\approx0.1414$ & $\frac1{32}\approx0.0313$\\[3pt]
ACIG & $\frac{2\sqrt2}{3}\bigl(1-s^{3/2}-(1-s)^{3/2}\bigr)$ & $\frac{2\sqrt2-2}{3}\approx0.2761$ & $\frac{2\sqrt2}{15}\approx0.1886$ & $\frac{\pi}{24}-\frac{4}{45}\approx0.0420$\\
\bottomrule
\end{tabular}
\end{table}

\section{Where does the fine structure live?}\label{sec:structure}

In all smooth examples of Section~\ref{sec:prelim} the tail order function is a product of powers $w_1^{a}w_2^{\kappa-a}$, determined by the tail order and one further number. Non-product forms occurred only for the Marshall--Olkin and the generalized Marshall--Olkin copulas (a kink along the diagonal), for the ACIG upper tail and for mixtures. This section explains this observation with the representation of \citet{WadsworthTawn2013}; see \citet{NoldeWadsworth2022} for its relation to other representations of multivariate extremes and to limit sets.

\subsection{The exponential scale}

For a copula $C$ with $C>0$ on $(0,1]^2$ define
\[
H_C(s)=-\log C(e^{-s_1},e^{-s_2}),\qquad s\in(0,\infty)^2 .
\]
Note that $C(e^{-s_1},e^{-s_2})=P(E_1\ge s_1,E_2\ge s_2)$ for $E_i=-\log U_i\sim\mathrm{Exp}(1)$, where $(U_1,U_2)\sim C$. Thus $H_C$ is the negative logarithm of a joint survival function on exponential margins, the object studied by \citet{WadsworthTawn2013}. Their basic assumption reads
\begin{equation}\label{eq:WT}
A_C(s):=\lim_{L\to\infty}\frac{H_C(Ls)}{L}\quad\text{exists for all }s\in(0,\infty)^2 .
\end{equation}
$A_C$ is $1$-homogeneous and nondecreasing with $A_C(s)\ge\max(s_1,s_2)$, and $A_C(s)\le s_1+s_2$ under positive quadrant dependence. The \emph{angular dependence function} is $\alpha_C(\omega)=A_C(\omega,1-\omega)$, $\omega\in(0,1)$ \citep[Section~2.2]{WadsworthTawn2013}, where it is denoted by $\lambda$. If $C\in\T$ satisfies \eqref{eq:WT}, then $C(u,u)\sim u^{\kappa}\ell(u)$ gives $H_C(L\one)=\kappa L-\log\ell(e^{-L})+o(1)$ as $L\to\infty$ (the term $o(1)$ is absent for the choice $\ell(u):=C(u,u)/u^{\kappa}$), and hence $\kappa_C=A_C(\one)=2\alpha_C(1/2)$, because $\log\ell(e^{-L})/L\to0$ for slowly varying $\ell$. This is the relation $1/\eta=2\alpha_C(1/2)$ of \citet{WadsworthTawn2013}.

\begin{lemma}\label{lem:logscale}
Let $C\in\T$. Then for every $x\in\R^2$
\[
g_C(x):=\lim_{L\to\infty}\bigl[H_C(L\one-x)-H_C(L\one)\bigr]\ \text{exists, and}\ b_C(e^{x_1},e^{x_2})=e^{-g_C(x)} .
\]
\end{lemma}

\begin{proof}
Put $u=e^{-L}$ and $w=(e^{x_1},e^{x_2})$, so that $uw_i=e^{-(L-x_i)}$. Then $C(uw)/C(u,u)=\exp\{-[H_C(L\one-x)-H_C(L\one)]\}$, and the left-hand side tends to $b_C(w)$ because $C(u,u)\sim u^{\kappa}\ell(u)$.
\end{proof}

Thus $b_C$ is the exponential of the negative \emph{local increment} of $H_C$ near the diagonal, while $A_C$ describes $H_C$ on the scale $L$. The question is when the first-order behaviour of $A_C$ at $\one$ already determines the local increments.

For a function $A$ on $(0,\infty)^2$ let $A'(\one;d)=\lim_{\eps\downarrow0}[A(\one+\eps d)-A(\one)]/\eps$ denote the one-sided directional derivative, whenever it exists, and define
\[
b_A(w)=\exp\{-A'(\one;-\log w)\},\qquad b^A(w)=\exp\{A'(\one;\log w)\},\qquad w\in(0,\infty)^2,
\]
with $\log w=(\log w_1,\log w_2)$.

\begin{lemma}\label{lem:bA}
Let $A$ be $1$-homogeneous, let $\kappa=A(\one)$, and suppose that $\alpha(\omega)=A(\omega,1-\omega)$ has one-sided derivatives $\alpha'_\pm$ at $\omega=1/2$. Put $a_\pm=(\kappa+\alpha'_\pm)/2$. Then
\[
b_A(w)=\begin{cases}w_1^{a_+}w_2^{\kappa-a_+},&w_1\le w_2,\\ w_1^{a_-}w_2^{\kappa-a_-},&w_1\ge w_2,\end{cases}\qquad
b^A(w)=\begin{cases}w_1^{a_-}w_2^{\kappa-a_-},&w_1\le w_2,\\ w_1^{a_+}w_2^{\kappa-a_+},&w_1\ge w_2.\end{cases}
\]
In particular $b_A=b^A=w_1^{a}w_2^{\kappa-a}$ with $a=\kappa/2+\alpha'(1/2)/2$ if $\alpha$ is differentiable at $1/2$.
\end{lemma}

\begin{proof}
For $d\in\R^2$ and small $\eps>0$ we have $A(\one+\eps d)=(2+\eps(d_1+d_2))\alpha(\omega_\eps)$ with $\omega_\eps=\tfrac12+\tfrac{\eps}{4}(d_1-d_2)+O(\eps^2)$. Hence
\[
A'(\one;d)=\frac{\kappa}{2}(d_1+d_2)+\frac{\alpha'_{\sigma}}{2}(d_1-d_2),\qquad \sigma=\begin{cases}+,&d_1>d_2,\\-,&d_1<d_2,\end{cases}
\]
where the second term vanishes for $d_1=d_2$. Inserting $d=-\log w$ resp.\ $d=\log w$ gives the claim.
\end{proof}

\subsection{Ray-regular tails}

\begin{theorem}\label{thm:rayregular}
Let $C$ satisfy \eqref{eq:WT}, and let $\alpha_C$ have one-sided derivatives $\alpha'_\pm$ at $1/2$. Assume that $R_C:=H_C-A_C$ satisfies the \emph{ray-regularity condition}
\begin{equation}\label{eq:RR}
\lim_{L\to\infty}\bigl[R_C(L\one-x)-R_C(L\one)\bigr]=0\qquad\text{for all }x\in\R^2 .
\end{equation}
Then $C\in\T$ with $\kappa_C=2\alpha_C(1/2)$, $\ell_C(u)=\exp\{-R_C(\log(1/u)\one)\}$ and $b_C=b_{A_C}$. In particular, $b_C$ is a product of powers if and only if $\alpha_C$ is differentiable at $1/2$; otherwise it is a product broken along the diagonal, with exponent $a_+$ of $w_1$ on $\{w_1\le w_2\}$ and $a_-$ on $\{w_1\ge w_2\}$.
\end{theorem}

\begin{proof}
By homogeneity $A_C(L\one-x)=LA_C(\one-x/L)$, so
\[
H_C(L\one-x)-H_C(L\one)=L\bigl[A_C(\one-x/L)-A_C(\one)\bigr]+\bigl[R_C(L\one-x)-R_C(L\one)\bigr]\longrightarrow A_C'(\one;-x),
\]
where the directional derivative exists by Lemma~\ref{lem:bA}. Moreover $C(u,u)=e^{-H_C(L\one)}=u^{\kappa}e^{-R_C(L\one)}$ with $L=\log(1/u)$, and $\ell(u)=e^{-R_C(L\one)}$ is slowly varying by \eqref{eq:RR} with $x=(\log t)\one$. The argument of Lemma~\ref{lem:logscale} gives $C(uw)/(u^{\kappa}\ell(u))\to\exp\{-A'_C(\one;-\log w)\}=b_{A_C}(w)$. The final statement follows from Lemma~\ref{lem:bA}, since $a_+\ne a_-$ unless $\alpha'_+=\alpha'_-$.
\end{proof}

\begin{remark}\label{rem:WT}
Condition \eqref{eq:RR} is the condition (3.2) of \citet{WadsworthTawn2013} at the ray $\omega=1/2$, combined with slow variation of $\ell$ along the diagonal. They describe it as a smoothness condition analogous to ray independence in \citet{LedfordTawn1997}. For differentiable $\alpha$, their conditional limit (3.3) at $\omega=1/2$ and the identification $c_1=\alpha(1/2)+\alpha'(1/2)/2$ in their Section~4.1 contain the product case of Theorem~\ref{thm:rayregular}, phrased for exceedances. In the representation $P(X_P>tx,Y_P>ty)\sim\mathcal L(tx,ty)\,t^{-1/\eta}x^{-c_1}y^{-c_2}$ of \citet{LedfordTawn1997}, the tail order function is $b(w)\propto w_1^{c_1}w_2^{c_2}g(1/w_1,1/w_2)$ with the ray dependence function $g$. A product form of $b$ is therefore the same as ray independence; models without ray independence are studied by \citet{RamosLedford2009}. What is new in Theorem~\ref{thm:rayregular} is the formulation for the tail order function on all of $(0,\infty)^2$ and, more importantly, the case of a kinked angular dependence function. \citet{WadsworthTawn2013} discuss non-differentiability only for asymptotically dependent vectors. The kinked case also occurs under asymptotic independence: for extreme-value copulas $H_C=A_C$ exactly, so \eqref{eq:RR} holds trivially, and for the Marshall--Olkin copula $\alpha(\omega)=1-\beta\min(\omega,1-\omega)$, $\alpha'_\pm=\pm\beta$, which reproduces Example~\ref{ex:families}(1). The same holds for the generalized Marshall--Olkin copula $C^{\mathrm{GMO}}_\kappa$ of Example~\ref{ex:families}(4): the explicit formula for $P(X_1>x_1,X_2>x_2)$ there gives $A_C(s)=\min(\kappa\max(s_1,s_2),s_1+s_2)$, which equals $\kappa\max(s_1,s_2)$ near $\one$ because $\kappa<2$, so $\alpha'_\pm=\pm\kappa$, $a_+=\kappa$, $a_-=0$; moreover $R_C(s)=H_C(s)-A_C(s)\to0$ as $\min(s)\to\infty$ with $|s_1-s_2|$ bounded, so \eqref{eq:RR} holds; Theorem~\ref{thm:rayregular} gives $b_C=b_{A_C}=\min(w_1,w_2)^{\kappa}$, in accordance with Example~\ref{ex:families}(4).
\end{remark}

\subsection{Log-concave tails}

Condition \eqref{eq:RR} concerns the remainder $R_C$ and may be hard to verify. The next result replaces it by a shape condition.

\begin{theorem}\label{thm:logconcave}
Let $C\in\T$ with $C>0$ on $(0,1]^2$ satisfy \eqref{eq:WT}, and assume that
\[
s\longmapsto C(e^{-s_1},e^{-s_2})=P(E_1\ge s_1,E_2\ge s_2)\quad\text{is log-concave on }(0,\infty)^2,
\]
i.e.\ $H_C$ is convex. Then $A_C$ and $\alpha_C$ are convex, $\alpha'_-\le\alpha'_+$, and
\[
b_{A_C}\le b_C\le b^{A_C}\qquad\text{on }(0,\infty)^2 .
\]
In particular, if $\alpha_C$ is differentiable at $1/2$, then $b_C(w)=w_1^{a}w_2^{\kappa-a}$ with $a=\kappa/2+\alpha_C'(1/2)/2$.
\end{theorem}

\begin{proof}
$A_C$ is a pointwise limit of the convex functions $s\mapsto H_C(Ls)/L$, hence convex; so is its restriction $\alpha_C$, whence the one-sided derivatives exist and $\alpha'_-\le\alpha'_+$. Fix $x\in\R^2$ and $\eps\in(0,1/\|x\|_\infty)$ (any $\eps>0$ if $x=0$), and let $L>1/\eps$. The function $h(t)=H_C(L\one-tx)$ is convex on $[-\eps L,\eps L]$, and the slopes of its chords increase. Therefore
\[
\frac{h(0)-h(-\eps L)}{\eps L}\ \le\ h(1)-h(0)\ \le\ \frac{h(\eps L)-h(0)}{\eps L}.
\]
As $L\to\infty$, the middle term tends to $g_C(x)$ by Lemma~\ref{lem:logscale}, while by \eqref{eq:WT} and homogeneity the outer terms tend to $[A_C(\one)-A_C(\one+\eps x)]/\eps$ and $[A_C(\one-\eps x)-A_C(\one)]/\eps$. Letting $\eps\downarrow0$ yields
\[
-A_C'(\one;x)\le g_C(x)\le A_C'(\one;-x),
\]
which is the asserted sandwich after exponentiation. If $\alpha_C$ is differentiable at $1/2$, then $A'_C(\one;\cdot)$ is linear by Lemma~\ref{lem:bA} and both bounds coincide.
\end{proof}

\begin{remark}\label{rem:prekopa}
(a) By Pr\'ekopa's theorem \citep{Prekopa1973}, the hypothesis holds whenever $(E_1,E_2)$ has a log-concave density, i.e.\ whenever $s\mapsto\log c(e^{-s_1},e^{-s_2})-s_1-s_2$ is concave, $c$ being the copula density; indeed $P(E\ge s)$ is the measure of the translate $s+[0,\infty)^2$ of a convex set. For extreme-value copulas $H_C=A_C$ is convex by definition.

(b) In contrast to Theorem~\ref{thm:rayregular}, Theorem~\ref{thm:logconcave} does not assume anything about the remainder $R_C$. It shows that, under log-concavity, smoothness of the angular dependence function at the diagonal alone forces ray independence.
\end{remark}

\subsection{Consequences for the order and the measures}

\begin{corollary}\label{cor:regularorder}
Let $C_1,C_2\in\T$ have the same tail order, $\ell_1/\ell_2\to\theta\in(0,\infty)$, and $b_{C_i}=b_{A_{C_i}}$ (as in Theorem~\ref{thm:rayregular}, or Theorem~\ref{thm:logconcave} with differentiable $\alpha_{C_i}$). Denote the one-sided derivatives of $\alpha_{C_i}$ at $1/2$ by $\alpha'_{i,\pm}$. Then
\[
C_1\preceq C_2\iff\theta\le1,\quad\alpha'_{1,+}\ge\alpha'_{2,+}\quad\text{and}\quad\alpha'_{1,-}\le\alpha'_{2,-}.
\]
In particular, if both angular dependence functions are differentiable at $1/2$, then $C_1$ and $C_2$ are comparable if and only if $\alpha_{C_1}'(1/2)=\alpha_{C_2}'(1/2)$, and then they are ordered by $\theta$.
\end{corollary}

\begin{proof}
By Proposition~\ref{thm:levels}, $C_1\preceq C_2$ iff $\theta b_1\le b_2$. For $w=(r,1)$ with $r\le1$, Lemma~\ref{lem:bA} gives $\theta b_1/b_2=\theta r^{(\alpha'_{1,+}-\alpha'_{2,+})/2}$, which is bounded by $1$ on $(0,1]$ iff $\theta\le1$ and $\alpha'_{1,+}\ge\alpha'_{2,+}$. The case $w=(1,r)$ gives $\theta r^{-(\alpha'_{1,-}-\alpha'_{2,-})/2}$, and homogeneity reduces every $w$ to one of these two cases.
\end{proof}

Within the ray-regular class, the finest level of the order is thus finite-dimensional and given by the triple $(\theta,\alpha'_+,\alpha'_-)$. A larger right derivative $\alpha'_+$ and a smaller left derivative $\alpha'_-$, i.e.\ a sharper kink of $\alpha$ at the diagonal, mean less tail dependence off the diagonal. In Example~\ref{ex:chain}, $\alpha'_+$ equals $0$, $1/2$ and $3/2$ for the Gumbel, the Marshall--Olkin and the generalized Marshall--Olkin copula, respectively (for the latter, $A(s)=\min(\kappa\max(s_1,s_2),s_1+s_2)$). The level measures also become explicit.

\begin{corollary}\label{cor:closedform}
Under the assumptions of Corollary~\ref{cor:regularorder}, with tail profile $\tau=\theta\,b_{A}$ and $a_\pm=(\kappa+\alpha'_\pm)/2\in[0,\kappa]$,
\[
\mu_1=\theta\Bigl[B_{1/2}(a_++1,\kappa-a_++1)+B_{1/2}(\kappa-a_-+1,a_-+1)\Bigr],
\]
where $B_x(p,q)=\int_0^xs^{p-1}(1-s)^{q-1}ds$ is the incomplete beta function. If $a_+\ge\kappa/2\ge a_-$ (e.g.\ if $\alpha$ is convex and symmetric), then $\mu_\infty=t(1/2)=\theta2^{-\kappa}$.
\end{corollary}

\begin{proof}
For $s\le1/2$ the profile is $t(s)=\theta s^{a_+}(1-s)^{\kappa-a_+}$, and for $s\ge1/2$ it is $\theta s^{a_-}(1-s)^{\kappa-a_-}$; integrate (substituting $s\mapsto1-s$ in the second part). The function $s^{a}(1-s)^{\kappa-a}$ increases on $[0,a/\kappa]$, which contains $[0,1/2]$ if $a\ge\kappa/2$; symmetrically on $[1/2,1]$.
\end{proof}

\subsection{Beyond ray regularity}

For ray-regular tails, and for log-concave tails with $\alpha$ differentiable at $1/2$, the tail order function is determined by $\kappa$ and $\alpha'_\pm$. Outside these classes it can carry additional information. Two mechanisms produce this.

\begin{example}[Mixtures]\label{ex:mixture}
Let $C=pC^{\mathrm{Gu}}_\kappa+(1-p)C^{\mathrm{GMO}}_\kappa$ with $p\in(0,1)$ and $\kappa\in(1,2)$. Since $H_C$ is a log-sum-exp, $A_C=\min(A_{\mathrm{Gu}},A_{\mathrm{GMO}})=A_{\mathrm{Gu}}$, because $(s_1^{\delta}+s_2^{\delta})^{1/\delta}\le\min(\kappa\max(s),s_1+s_2)$ for $2^{1/\delta}=\kappa$. Hence $\alpha_C$ is differentiable at $1/2$ with $\alpha_C'(1/2)=0$ and $b_{A_C}=(w_1w_2)^{\kappa/2}$. However, by Example~\ref{ex:families}(6),
\[
b_C(w)=p(w_1w_2)^{\kappa/2}+(1-p)\min(w_1,w_2)^{\kappa}<b_{A_C}(w)\qquad(w_1\ne w_2).
\]
Thus $b_C$ is not of product form and violates the lower bound of Theorem~\ref{thm:logconcave}; consequently $C$ is neither ray-regular nor log-concave on the exponential scale. The angular dependence function only records the dominating component, while $b_C$ records the mixture.
\end{example}

\begin{example}[Hidden regular variation: the ACIG upper tail]\label{ex:acig}
For $C^{\mathrm{ig}}_\kappa$, $\kappa\in(1,2)$, Lemma~\ref{lem:acig}(3) gives
\[
\alpha(\omega)=\max(\omega,1-\omega)+(\kappa-1)\min(\omega,1-\omega),
\]
so $\alpha$ is convex and kinked at $1/2$ with $\alpha'_\pm=\pm(2-\kappa)$. Lemma~\ref{lem:bA} gives $b_A(r,1)=r$ and $b^A(r,1)=r^{\kappa-1}$ for $r\le1$, whereas
\[
b(r,1)=\frac{(1+r)^{\kappa}-1-r^{\kappa}}{2^{\kappa}-2}
\]
is smooth and not a broken product. Hence \eqref{eq:RR} fails. Nevertheless the sandwich $b_A\le b\le b^A$ of Theorem~\ref{thm:logconcave} holds for every $\kappa\in(1,2)$ (Lemma~\ref{lem:acigineq}(1)); we do not know whether $C^{\mathrm{ig}}_\kappa$ is log-concave on the exponential scale. Here the non-product shape comes from the regularly varying frailty, i.e.\ from hidden regular variation in the sense of \citet{HuaJoeLi2014}.
\end{example}

\section{Discussion}\label{sec:discussion}

We extended the program of \citet{SiburgStrothmannWeiss2024} from the tail dependence function to the tail expansion \eqref{eq:expansion}, which remains informative under asymptotic independence. The asymptotic tail order $\preceq$ refines $\le_{td}$ and is lexicographic in tail order, slowly varying part and tail order function. This structure forces a two-tier notion of measures: global measures that are functions of the tail order, and level measures that compare copulas with the same tail order and comparable slowly varying parts. For ray-regular tails, and for log-concave tails with an angular dependence function that is differentiable at the diagonal, the finest level is described by the one-sided slopes of the angular dependence function at the diagonal; for general log-concave tails these slopes still bound the tail order function from both sides. A genuinely functional tail order function, beyond (broken) products, appears for mixtures and hidden regular variation.

Several questions remain open.
\begin{enumerate}
\item \emph{Construction via extreme-value copulas.} For $\kappa=1$, \citet[Theorem~3.4]{SiburgStrothmannWeiss2024} turn every concordance measure into a tail dependence measure via the extreme-value copula $C^{EV}_\Lambda$. An analogue for $\kappa>1$, mapping a tail profile to a copula in a monotone way, is not known to us.
\item \emph{Characterization of tail profiles.} By \citet[Theorem~8.2.3]{Jaworski2010}, locally Lipschitz, grounded, $2$-increasing homogeneous functions of degree $\kappa>1$ are tail profiles, while leading parts of tail expansions of degree $\kappa>1$ need not be continuous \citep[Corollary~8.3.1]{Jaworski2010}. A characterization of all tail profiles, and of the admissible pairs $(\alpha,b)$, would sharpen Sections~\ref{sec:measures} and~\ref{sec:structure}.
\item \emph{Local dominance.} Under which conditions does $C_1\prec C_2$ imply $C_1\le_{loc}C_2$, in analogy to \citet{SiburgStrothmann2024}? Proposition~\ref{prop:cone} gives dominance on cones; the behaviour near the axes is governed by other directions in the representation of \citet{WadsworthTawn2013}.
\item \emph{Higher dimensions.} The definitions extend verbatim to $d\ge2$. The structure results of Section~\ref{sec:structure} require the $d$-dimensional angular dependence function and the subcone structure of hidden regular variation \citep{HuaJoeLi2014}.
\item \emph{Statistics.} The level measures depend on $b$, whose plug-in estimator $\widehat C(uw)/\widehat C(u,u)$ converges slowly. For the Gaussian copula, the proof of Lemma~\ref{lem:gauss} shows that the relative error of $C(uw)/C(u,u)$ as an approximation of $b(w)$ is of order $1/\log(1/u)$ for fixed $w$; numerically it is still about $1\%$ at $u=10^{-48}$ for $\rho=1/2$ and $w=(0.1,1)$. Simulation studies and data analyses are left for future work.
\end{enumerate}

\appendix
\section{Proofs for the examples}\label{app:proofs}

\begin{lemma}[Gaussian copula]\label{lem:gauss}
Let $C$ be the Gaussian copula with correlation $\rho\in(-1,1)$ and $\kappa=2/(1+\rho)$. Then, for every $w\in(0,\infty)^2$,
\[
C(uw_1,uw_2)\sim c_\rho\,\bigl(2\log(1/u)\bigr)^{-\rho/(1+\rho)}\,u^{\kappa}\,(w_1w_2)^{\kappa/2},\qquad u\downarrow0,
\]
with $c_\rho$ as in \eqref{eq:gaussell}. Moreover $C(uw)/C(u,u)=(w_1w_2)^{\kappa/2}\bigl(1+O(1/\log(1/u))\bigr)$ for fixed $w$.
\end{lemma}

\begin{proof}
Let $(X,Y)$ be standard bivariate normal with correlation $\rho$ and density
\[
\varphi_2(x,y)=\frac{1}{2\pi\sqrt{1-\rho^2}}\exp\Bigl\{-\frac{q(x,y)}{2(1-\rho^2)}\Bigr\},\qquad q(x,y)=x^2-2\rho xy+y^2 .
\]

\emph{Step 1 (orthant probabilities).} For $a_1,a_2>0$ put $\gamma_1=(a_1-\rho a_2)/(1-\rho^2)$ and $\gamma_2=(a_2-\rho a_1)/(1-\rho^2)$. Since $q(a_1+s,a_2+t)=q(a_1,a_2)+2(1-\rho^2)(\gamma_1 s+\gamma_2 t)+q(s,t)$,
\[
P(X>a_1,Y>a_2)=\varphi_2(a_1,a_2)\int_0^\infty\!\!\int_0^\infty e^{-\gamma_1 s-\gamma_2 t-Q(s,t)/2}\,ds\,dt,\qquad Q=\frac{q}{1-\rho^2}\ge0 .
\]
Substituting $s=\sigma/\gamma_1$, $t=\tau/\gamma_2$ and using $0\le1-e^{-Q/2}\le Q/2$, we get for $\gamma_1,\gamma_2>0$ (cf.\ the bounds of \citet{Savage1962})
\begin{equation}\label{eq:orthant}
P(X>a_1,Y>a_2)=\frac{\varphi_2(a_1,a_2)}{\gamma_1\gamma_2}\bigl(1+O(\gamma_1^{-2}+\gamma_2^{-2})\bigr).
\end{equation}

\emph{Step 2 (quantiles).} Let $a=a(u)=-\Phi^{-1}(u)$. The Mills ratio expansion $\Phi(-a)=\varphi(a)a^{-1}(1+O(a^{-2}))$ gives
\begin{equation}\label{eq:quantile}
\log u=-\tfrac12a^2-\log a-\tfrac12\log(2\pi)+O(a^{-2}),
\end{equation}
hence $a^2\sim2\log(1/u)$. Let $a_i=a(uw_i)$. Subtracting \eqref{eq:quantile} at $u$ and at $uw_i$ gives $a_i^2-a^2=-2\log w_i-2\log(a_i/a)+O(a^{-2})$. Since $a_i^2\sim2\log(1/u)\sim a^2$, we have $a_i/a\to1$; hence $a_i^2-a^2$ is bounded, $\log(a_i/a)=\frac12\log(a_i^2/a^2)=O(a^{-2})$, and $a_i^2-a^2=-2\log w_i+O(a^{-2})$. Consequently $a(a_i-a)=(a_i^2-a^2)\,a/(a+a_i)=-\log w_i+O(a^{-2})$ and $(a_1-a)(a_2-a)=O(a^{-2})$.

\emph{Step 3.} Since $(X,Y)\overset{d}{=}(-X,-Y)$, $C(uw)=P(X>a_1,Y>a_2)$. As $\gamma_1,\gamma_2\sim a/(1+\rho)\to\infty$, \eqref{eq:orthant} applies to $C(uw)$ and $C(u,u)$, and the prefactors $\gamma_1\gamma_2$ have ratio $1+O(a^{-2})$. Moreover
\[
\begin{aligned}
q(a_1,a_2)-q(a,a)&=\sum_i(a_i^2-a^2)-2\rho\bigl[a(a_1-a)+a(a_2-a)+(a_1-a)(a_2-a)\bigr]\\
&=-2(1-\rho)\log(w_1w_2)+O(a^{-2}),
\end{aligned}
\]
so
\[
\frac{C(uw)}{C(u,u)}=\exp\Bigl\{-\frac{q(a_1,a_2)-q(a,a)}{2(1-\rho^2)}\Bigr\}\bigl(1+O(a^{-2})\bigr)=(w_1w_2)^{1/(1+\rho)}\bigl(1+O(a^{-2})\bigr).
\]
On the diagonal, $\gamma_1\gamma_2=a^2/(1+\rho)^2$ and $q(a,a)=2(1-\rho)a^2$, so
\[
C(u,u)\sim\frac{(1+\rho)^2}{2\pi\sqrt{1-\rho^2}\,a^2}\,e^{-a^2/(1+\rho)}.
\]
By \eqref{eq:quantile}, $e^{-a^2/2}=\sqrt{2\pi}\,a\,u\,(1+O(a^{-2}))$, hence $e^{-a^2/(1+\rho)}=(e^{-a^2/2})^{\kappa}\sim(2\pi)^{\kappa/2}a^{\kappa}u^{\kappa}$ and
\[
C(u,u)\sim(1+\rho)^2(1-\rho^2)^{-1/2}(2\pi)^{\kappa/2-1}\,a^{\kappa-2}\,u^{\kappa}.
\]
Finally $a^{\kappa-2}\sim(2\log(1/u))^{(\kappa-2)/2}$ and $(\kappa-2)/2=-\rho/(1+\rho)$. The error terms in the ratio $C(uw)/C(u,u)$ are $O(a^{-2})=O(1/\log(1/u))$. In the last step, \eqref{eq:quantile} gives
\[
a^2=2\log(1/u)\Bigl(1+O\Bigl(\frac{\log\log(1/u)}{\log(1/u)}\Bigr)\Bigr).
\]
\end{proof}

\begin{lemma}[ACIG copula]\label{lem:acig}
Let $\kappa\in(1,2)$, $\psi(s)=2\Gamma(\kappa)^{-1}s^{\kappa/2}K_\kappa(2\sqrt s)$ the Laplace transform of $1/G$ with $G\sim\mathrm{Gamma}(\kappa,1)$, $C_\psi(v_1,v_2)=\psi(\psi^{-1}(v_1)+\psi^{-1}(v_2))$, and $C^{\mathrm{ig}}_\kappa(u_1,u_2)=u_1+u_2-1+C_\psi(1-u_1,1-u_2)$. Put $D(x,y)=(x+y)^{\kappa}-x^{\kappa}-y^{\kappa}$.
\begin{enumerate}
\item $h(s):=\psi(s)-1+s/(\kappa-1)=c_hs^{\kappa}+r(s)$ with $c_h=-\Gamma(1-\kappa)/\Gamma(1+\kappa)>0$, and there is $M$ with $|r(s_1+s_2)-r(s_1)-r(s_2)|\le Ms_1s_2$ for $s_1+s_2\le1$.
\item $C^{\mathrm{ig}}_\kappa\in\T$ with tail order $\kappa$, $\ell\to\Upsilon_L(C^{\mathrm{ig}}_\kappa)=c_h(2^{\kappa}-2)(\kappa-1)^{\kappa}$ and $b=D/(2^{\kappa}-2)$.
\item $C^{\mathrm{ig}}_\kappa$ satisfies \eqref{eq:WT} with $\alpha(\omega)=\max(\omega,1-\omega)+(\kappa-1)\min(\omega,1-\omega)$.
\end{enumerate}
\end{lemma}

\begin{proof}
(1) For non-integer $\kappa$, $K_\kappa=\pi(I_{-\kappa}-I_\kappa)/(2\sin(\pi\kappa))$ and $I_\nu(2\sqrt s)=\sum_{k\ge0}s^{k+\nu/2}/(k!\,\Gamma(k+\nu+1))$. With $\pi/(\Gamma(\kappa)\sin(\pi\kappa))=\Gamma(1-\kappa)$ this gives
\[
\psi(s)=\Gamma(1-\kappa)\Bigl[\sum_{k\ge0}\frac{s^{k}}{k!\,\Gamma(k+1-\kappa)}-\sum_{k\ge0}\frac{s^{k+\kappa}}{k!\,\Gamma(k+1+\kappa)}\Bigr].
\]
The coefficients of $s^0$, $s^1$ and $s^{\kappa}$ are $1$, $\Gamma(1-\kappa)/\Gamma(2-\kappa)=-1/(\kappa-1)$ and $c_h$; $c_h>0$ because $\Gamma(1-\kappa)<0$ for $\kappa\in(1,2)$. Thus $r(s)=\sum_{k\ge2}a_ks^k+\sum_{k\ge1}b_ks^{k+\kappa}$ with absolutely convergent series. For $p(s)=s^m$, $m\ge2$, we have $p(s_1+s_2)-p(s_1)-p(s_2)=\int_0^{s_1}\!\!\int_0^{s_2}p''(x+y)\,dy\,dx\le m(m-1)s_1s_2$ if $s_1+s_2\le1$; summing gives the bound with $M=\sum_k|a_k|k(k-1)+\sum_k|b_k|(k+\kappa)(k+\kappa-1)<\infty$.

(2) Let $s_i=\psi^{-1}(1-u_i)$. Then $u_i=1-\psi(s_i)=s_i/(\kappa-1)-h(s_i)$ and
\[
C^{\mathrm{ig}}_\kappa(u_1,u_2)=u_1+u_2-1+\psi(s_1+s_2)=h(s_1+s_2)-h(s_1)-h(s_2)=c_hD(s_1,s_2)+O(s_1s_2).
\]
From $u_i=s_i/(\kappa-1)-c_hs_i^{\kappa}+O(s_i^2)$ we get $s_i=(\kappa-1)u_i(1+O(u_i^{\kappa-1}))$. For $u_i=uw_i$, homogeneity and continuity of $D$ give $c_hD(s_1,s_2)=c_h(\kappa-1)^{\kappa}u^{\kappa}D(w)(1+o(1))$, while $O(s_1s_2)=O(u^2)=o(u^{\kappa})$. Hence $C^{\mathrm{ig}}_\kappa(uw)\sim c_h(\kappa-1)^{\kappa}D(w)\,u^{\kappa}$, which is the claim, since $D(\one)=2^{\kappa}-2$.

(3) Writing $D(x,y)=\kappa(\kappa-1)\int_0^x\!\!\int_0^y(p+q)^{\kappa-2}\,dq\,dp$ and using $(2y)^{\kappa-2}\le(p+q)^{\kappa-2}$ and $\int_0^y(p+q)^{\kappa-2}dq\le y^{\kappa-1}/(\kappa-1)$, we obtain
\begin{equation}\label{eq:Dbounds}
\kappa(\kappa-1)2^{\kappa-2}\,xy^{\kappa-1}\le D(x,y)\le\kappa\,xy^{\kappa-1},\qquad0<x\le y .
\end{equation}
With $s_{\min}\le s_{\max}$ the remainder in (2) is at most $Ms_{\min}s_{\max}=s_{\min}s_{\max}^{\kappa-1}\cdot Ms_{\max}^{2-\kappa}$, so by \eqref{eq:Dbounds} $C^{\mathrm{ig}}_\kappa(u_1,u_2)=c_hD(s_1,s_2)(1+o(1))$ uniformly as $\max(u_1,u_2)\to0$. Put $u_i=e^{-t_i}$, so that $s_i=(\kappa-1)e^{-t_i}(1+o(1))$. Then \eqref{eq:Dbounds} yields
\[
H(t)=-\log C^{\mathrm{ig}}_\kappa(e^{-t_1},e^{-t_2})=\max(t_1,t_2)+(\kappa-1)\min(t_1,t_2)+O(1)
\]
uniformly as $\min(t_1,t_2)\to\infty$. Dividing $H(L\omega,L(1-\omega))$ by $L$ gives the claim.
\end{proof}

\begin{lemma}\label{lem:acigineq}
Let $\kappa\in(1,2)$ and $b$ be the tail order function of $C^{\mathrm{ig}}_\kappa$.
\begin{enumerate}
\item With $A$ and $\alpha$ from Lemma~\ref{lem:acig}(3), $b_A\le b\le b^A$; equivalently, $(2^{\kappa}-2)\,r\le D(r,1)\le(2^{\kappa}-2)\,r^{\kappa-1}$ for $r\in(0,1]$.
\item For $\kappa=3/2$: $\Upsilon_L(C^{\mathrm{ig}}_{3/2})\,b(w)\le\frac{7\sqrt2}{12}(w_1w_2)^{3/4}$ for all $w\in(0,\infty)^2$. In particular $C^{\mathrm{ig}}_{3/2}\prec C^{\mathrm{Gu}}$, where $C^{\mathrm{Gu}}$ is the Gumbel copula with tail order $3/2$.
\end{enumerate}
\end{lemma}

\begin{proof}
(1) Since $\alpha'_\pm=\pm(2-\kappa)$, Lemma~\ref{lem:bA} gives $a_+=1$ and $a_-=\kappa-1$, so $b_A(r,1)=r$ and $b^A(r,1)=r^{\kappa-1}$ for $r\le1$. By homogeneity and symmetry it suffices to consider $w=(r,1)$ with $r\in(0,1]$. Let $f(q)=D(q,1)/q$ for $q>0$. By the integral representation in the proof of Lemma~\ref{lem:acig}(3),
\[
f(q)=\frac1q\int_0^qG(p)\,dp,\qquad G(p)=\kappa(\kappa-1)\int_0^1(p+t)^{\kappa-2}\,dt,
\]
and $G$ is decreasing, so $f$ is nonincreasing as an average of a decreasing function. For $r\le1$ this gives $D(r,1)=rf(r)\ge rf(1)=(2^{\kappa}-2)r$. Moreover, by homogeneity and symmetry, $D(r,1)=r^{\kappa}D(1/r,1)=r^{\kappa-1}f(1/r)\le r^{\kappa-1}f(1)$, since $1/r\ge1$.

(2) For $\kappa=3/2$, $\Upsilon_L(C^{\mathrm{ig}}_{3/2})/(2^{3/2}-2)=c_h(1/2)^{3/2}=2\sqrt2/3$, so $\Upsilon_L(C^{\mathrm{ig}}_{3/2})\,b(r,1)=\frac{2\sqrt2}3D(r,1)$. By symmetry and homogeneity it suffices to show $\frac{2\sqrt2}3D(r,1)\le\frac{7\sqrt2}{12}r^{3/4}$ for $r\in(0,1]$. Since the third derivative of $(1+r)^{3/2}$ is negative, $(1+r)^{3/2}\le1+\frac32r+\frac38r^2$, hence $D(r,1)\le\frac32r+\frac38r^2-r^{3/2}$. With $x=r^{1/4}\in(0,1]$,
\[
\frac{D(r,1)}{r^{3/4}}\le p(x):=\tfrac32x+\tfrac38x^5-x^3 .
\]
Now $p'(x)=\frac{15}{8}\bigl(x^4-\frac85x^2+\frac45\bigr)>0$, because the quadratic $y^2-\frac85y+\frac45$ has negative discriminant. Thus $p(x)\le p(1)=\frac78$ and $\frac{2\sqrt2}3D(r,1)\le\frac{7\sqrt2}{12}r^{3/4}$. Since $\frac{7\sqrt2}{12}<1$ and both copulas have tail order $3/2$ with constant slowly varying parts, Proposition~\ref{thm:levels} gives $C^{\mathrm{ig}}_{3/2}\preceq C^{\mathrm{Gu}}$. The relation is strict because $\Upsilon_L(C^{\mathrm{ig}}_{3/2})\,b(\one)=(8-4\sqrt2)/3<1=b^{\mathrm{Gu}}(\one)$.
\end{proof}

\end{document}